\documentclass[12pt]{article}
\usepackage{amsmath,amsfonts}
\usepackage{verbatim}
\usepackage{latexsym}
\usepackage{graphicx}
\usepackage{color}
\makeatletter
\@addtoreset{equation}{section}

\begin{document}

\newcommand{\E}{\mathbb{E}}
\newcommand{\PP}{\mathbb{P}}
\newcommand{\RR}{\mathbb{R}}

\newtheorem{theorem}{Theorem}[section]
\newtheorem{lemma}[theorem]{Lemma}
\newtheorem{coro}[theorem]{Corollary}
\newtheorem{defn}[theorem]{Definition}
\newtheorem{assp}[theorem]{Assumption}
\newtheorem{expl}[theorem]{Example}
\newtheorem{prop}[theorem]{Proposition}
\newtheorem{rmk}[theorem]{Remark}

\newcommand\tq{{\scriptstyle{3\over 4 }\scriptstyle}}
\newcommand\qua{{\scriptstyle{1\over 4 }\scriptstyle}}
\newcommand\hf{{\textstyle{1\over 2 }\displaystyle}}
\newcommand\hhf{{\scriptstyle{1\over 2 }\scriptstyle}}

\newcommand{\proof}{\noindent {\it Proof}. }
\newcommand{\eproof}{\hfill $\Box$} % {\indent\vrule height6pt width4pt depth1pt\hfil\par\medbreak}

\def\a{\alpha} \def\g{\gamma}
\def\e{\varepsilon} \def\z{\zeta} \def\y{\eta} \def\o{\theta}
\def\vo{\vartheta} \def\k{\kappa} \def\l{\lambda} \def\m{\mu} \def\n{\nu}
\def\x{\xi}  \def\r{\rho} \def\s{\sigma}
\def\p{\phi} \def\f{\varphi}   \def\w{\omega}
\def\q{\surd} \def\i{\bot} \def\h{\forall} \def\j{\emptyset}

\def\be{\beta} \def\de{\delta} \def\up{\upsilon} \def\eq{\equiv}
\def\ve{\vee} \def\we{\wedge}

\def\F{{\cal F}}
\def\T{\tau} \def\G{\Gamma}  \def\D{\Delta} \def\O{\Theta} \def\L{\Lambda}
\def\X{\Xi} \def\S{\Sigma} \def\W{\Omega}
\def\M{\partial} \def\N{\nabla} \def\Ex{\exists} \def\K{\times}
\def\V{\bigvee} \def\U{\bigwedge}

\def\1{\oslash} \def\2{\oplus} \def\3{\otimes} \def\4{\ominus}
\def\5{\circ} \def\6{\odot} \def\7{\backslash} \def\8{\infty}
\def\9{\bigcap} \def\0{\bigcup} \def\+{\pm} \def\-{\mp}
\def\la{\langle} \def\ra{\rangle}

\def\tl{\tilde}
\def\trace{\hbox{\rm trace}}
\def\diag{\hbox{\rm diag}}
\def\for{\quad\hbox{for }}
\def\refer{\hangindent=0.3in\hangafter=1}

\newcommand\wD{\widehat{\D}}

\title{
\bf Advances in Stabilisation of Hybrid Stochastic Differential Equations by Delay Feedback Control
%\thanks{Partially supported by ...}
 }
\author{
{\bf Junhao Hu${}^{1}$,  Wei Liu${}^{2,}$\thanks{The corresponding author. Email: weiliu@shnu.edu.cn, lwbvb@hotmail.com}, Feiqi Deng${}^3$, Xuerong Mao${}^4$ }
\\
${}^1$ School of Mathematics and Statistics, \\
South-Central University for Nationalities, \\
 Wuhan, Hubei 430074, China.\\
${}^2$ Department of Mathematics, Shanghai Normal University, \\
Shanghai, 200234, China. \\
${}^3$ Systems Engineering Institute, \\
South China University of Technology, \\
Guangzhou 510640, China. \\
${}^4$ Department of Mathematics and Statistics, \\
University of Strathclyde, Glasgow G1 1XH, U.K. \\
 }

\date{}

\maketitle

\begin{abstract}
A novel approach to design the feedback control based on past states is proposed for hybrid stochastic differential equations (HSDEs). This new theorem builds up the connection between the delay feedback control and the control function without delay terms, which enables one to construct the delay feedback control using the existing results on stabilities of HSDEs. Methods to find the upper bound of the length of the time delay are also investigated. Numerical simulations are presented to demonstrate the new theorem.


\end{abstract}

\section{Introduction}

Many systems in the real word may experience abrupt changes
in their structures and parameters due to sudden changes of system factors, for example,  a failure of a power station in a network, a change of interest rate in an economic system, an environmental change in an ecological system.  Hybrid systems driven by continuous-time Markov chains have been
used widely to model these systems (see, e.g., \cite{And,BM04,Lew00,MMP,MY06,Mar}).

One important class of hybrid systems is the
hybrid stochastic differential equations
(SDEs; also known as SDEs with Markovian switching).
Hybrid SDEs are in general described by
\begin{equation} \label{1.1}
dx(t) = f(x(t),r(t),t)dt + g(x(t),r(t),t)dB(t).
\end{equation}
Here the state $x(t)$ takes values in $R^n$ and the mode $r(t)$ is described by a Markov chain taking values in a finite
space $S=\{1,2,\cdots, N\}$,  $B(t)$ is a Brownian motion,
 $f$ and $g$ are referred to as the drift and diffusion coefficient, respectively.
 (Further details on the notation will be stated in Section 2.)
 One of the important issues in the study of hybrid SDEs is the analysis of stability.
For example, Ji et al. \cite{JC}, and Mariton \cite{Mar}
 studied the stability of the jump linear systems.
Basak et al. \cite{BBG} discussed the stability of
semi-linear hybrid stochastic differential equations (SDEs)
while Mao \cite{M99} investigated the stability of a nonlinear
hybrid SDEs.
Shaikhet \cite{Sha} took the time delay into account and
considered the stability of
 semi-linear hybrid SDEs with delay, while   Mao et al. \cite{MMP} investigated the stability
of a nonlinear hybrid SDEs with delay.  Taking into account of the parameter uncertainty,
Mao \cite{M02} studied
the stability of hybrid stochastic delay interval systems.  There is now
an intensive literature in the area of hybrid SDEs (for further references
see, e.g., \cite{M99b,M07,MLXG,MYY,SMYI,SLXZ,WWSF,YH}).

 Given an $n$-dimensional unstable hybrid SDE in the form of (\ref{1.1}), it is classical to
find a feedback control $u(x(t),r(t),t)$, based on the current state $x(t)$,
for the controlled system
\begin{equation} \label{1.2}
dx(t) = [f(x(t),r(t),t)+ u(x(t),r(t),t) ]dt + g(x(t),r(t),t)dB(t)
\end{equation}
to become stable. This stabilisation problem by (non-delay) state feedback controls has been well studied (see, e.g., \cite{MYY,MY06,YM04}).

On the other hand, there is always a time lag between the time when the observation of the state is made and the time when the feedback control reaches the system. (It takes 1.28 seconds for a radio signal from the moon to reach the earth.)  Traditionally, we usually assume that the time lag is extremely small (namely the
feedback control acts extremely fast).  Nevertheless, there is a delay, say $\e$.  So the real controlled system should
be in the form of
\begin{equation} \label{1.2a}
dx(t) = [f(x(t),r(t),t)+ u(x(t-\e),r(t),t) ]dt + g(x(t),r(t),t)dB(t).
\end{equation}
For example, the time unit is of year while the time lag $\e$ is
%a second ($=3.17*10^{-8}$ year) or
%a millisecond ($=3.17*10^{-11}$ year) or
a microsecond ($=3.17*10^{-14}$ year)
or a nanosecond ($=3.17*10^{-17}$ year).
It seems that it has been taken as a grant that if the controlled system (\ref{1.2}) is stable, so is (\ref{1.2a}) provided $\e$ is sufficiently small.  However, there is a counter example which shows that
the controlled system (\ref{1.2}) is exponentially stable in $p$th moment, but the corresponding system (\ref{1.2a}) is NOT no matter how small $\e$ is (please see the Appendix).

More usefully, the time lag may NOT be extremely small.  In practice, to reduce the control cost, it may be more implementable  to allow the feedback control to act reasonably fast but not necessarily extremely fast.
Let $\T$ denote the time lag between the time when the observation of the state is made and the time when the feedback control reaches the system.  It is then certainly more
realistic that the control should depend on the past state $x(t-\T)$.  Accordingly, the control should be of
the form $u(x(t-\T), r(t),t)$.  Consequently,
 the stabilisation problem  becomes  to design a delay feedback
 control $u(x(t-\T),r(t),t)$  for the controlled system
\begin{equation} \label{1.3}
 dx(t)= [f(x(t),r(t),t)+ u(x(t-\T),r(t),t) ]dt + g(x(t),r(t),t)dB(t)
\end{equation}
to become stable.  Mao, Lam and Huang \cite{MLH} were the first to study this stabilisation problem by the delay feedback control for hybrid SDEs,
although the method of delay feedback controls has been well used in the area of ordinary differential equations (see, e.g., \cite{AU04,CLH05,Pyr95}).
The main idea in \cite{MLH} was to use the theory of linear matrix inequalities (LMIs) along with the method of Lyapunov functionals to design the \emph{linear} delay feedback control in the form of $u(x(t-\T),r(t),t) = F(r(t))G(r(t))x(t-\T)$, where $F:S\to R^{n\K l}$ and $G: S\to R^{l\K n}$. They discussed two cases: (i) state feedback, namely design $F(\cdot)$ when $G(\cdot)$ was given; (ii) output injection, namely design $G(\cdot)$ when $F(\cdot)$ was given. The advantage of the results in \cite{MLH} lies in that either $F(\cdot)$ or $G(\cdot)$ can be solved efficiently by the technique of LMIs. The disadvantage is that the given unstable hybrid SDE should be either linear or nonlinear but dominated by a linear system.
%(e.g., there are $N+1$ symmetric $n\K n$-matrices $Z$ and $Z_i$ ($i\in S$) with $Z$ being positive-definite  such that $2x^TZf(x,i,t)+g^T(x,i,t)Zg(x,i,t)\le z^Z_ix$)
Such a limitation is unavoidable due to the key technique of LMIs used.

In this paper we will take a completely new approach in order to avoid this limitation. As we mentioned before, the stabilization by the (non-delay) feedback control (namely, the controlled system (\ref{1.2}) has been well studied. In other words, there are lots of results on how to design the control function $u(x,i,t)$ to make the controlled system (ˆ\ref{1.2}) stable.

\begin{itemize}

\item[] \textbf{Question}: \emph{Can we make use of this same control function to
make the delay feedback controlled system (\ref{1.3}) stable? }

\end{itemize}

\noindent
Should the answer is yes, this would be great as the stabilization problem (\ref{1.3}) by a delay feedback control could be transferred to the well-known classical
stabilization problem (\ref{1.2}) by a  non-delay feedback control.  Is this possible?   To see the possibility, we rearrange the controlled system (\ref{1.3}) as
\begin{align} \label{1.4}
 dx(t) &= [f(x(t),r(t),t)+ u(x(t),r(t),t) ]dt + g(x(t),r(t),t)dB(t)
 \nonumber \\
 & -[u(x(t),r(t),t)-u(x(t-\T),r(t),t)]dt.
\end{align}
Comparing this with (\ref{1.2}), we may regard it (i.e., system (\ref{1.3})) as the perturbed system of (\ref{1.2}) with the perturbation $-[u(x(t),r(t),t)-u(x(t-\T),r(t),t)]dt$.  If the time lag $\tau$ (the  duration between the time when the state observation is made and the time when the feedback control reaches the system) is sufficiently small (namely, the feedback control acts sufficiently fast) while the control function $u(x,i,t)$ is globally Lipschitz continuous in $x$,  then the perturbation might be sufficiently small so that system (\ref{1.4}) should perform in a similar way as system (\ref{1.2}) does (namely, stable).  It is this perturbation idea that motivates us to show in this paper the following result:

\begin{itemize}

\item[] \textbf{Answer}: \emph{Under the global Lipschitz condition on the system coefficients $f, g$ and the control function $u$, if the control function $u$ makes the controlled system (\ref{1.2}) to be exponentially stable in the $p$th moment ($p>0$)
    then there is a positive number $\tau^*$, which can be determined numerically, such that the same control function $u$ will also make the controlled system (\ref{1.3}) to be exponentially stable in the $p$th moment as long as the delay feedback control $u(x(t-\tau),r(t),t)$ acts sufficiently fast in the sense $\tau < \tau^*$. }

\end{itemize}

Let us highlight a couple of important features of this new result:

\begin{itemize}

\item The result covers a much wider class of nonlinear hybrid SDEs than Mao et al. \cite{MLH}.

\item The delay feedback stabilisation problem (\ref{1.3}) is transferred to the classical stabilisation problem (\ref{1.2}) so that many existing results and techniques can be used to design the required control function $u(x,i,t)$.

\item The positive number $\tau^*$ can be determined numerically, which means our theory can be implemented easily.

\item More importantly, our new result gives a theoretical support for the general practice of the non-delay state feedback control. In practice,
    there is always a time lag, though it might be extremely small, from the moment when the state is observed
    to the time when the state feedback control reaches the system.
    In other words, the controlled system in practice should be in the form of (\ref{1.2a}).
     In the past, we have always designed the control function to make
     the controlled system (\ref{1.2}) to be stable while required the feedback control to act extremely quickly (like a non-delay though there is a delay), namely $\e$ is extremely small whence less than $\tau^*$
     (it seems that the physical existence of $\tau^*$ was not known before). Our new theory explains why this general practice has worked for so many years (under the global Lipschitz condition of course).

  \item Mao et al. \cite{MLH} only studied the stabilisation in the mean square exponential stability. We here discuss the stabilisation in more general $p$th exponential stability (for $p>0$).  This is particularly significant when $p\in (0,1)$ where
  the stabilisation effect of Brownian motions could be used (see Example \ref{E5.3} for more details please).

\end{itemize}

Let us begin to establish our new theory.  We will state some preliminaries in Section 2 and present a number of lemmas in Section 3.  We will prove the main results of this paper in Section 4 while demonstrate how our new results can be implemented easily in Section 5 making use of the existing theory on the non-delay state feedback controls.  We will conclude our paper finally in Section 6.

\section{Preliminaries}

Throughout this paper, unless otherwise specified, we will use the following notation.
If $A$ is a vector or matrix, its transpose is denoted by $A^T$.
If $x\in R^n$, then $|x|$ is its Euclidean norm.
If $A$ is a matrix, we let $|A| = \sqrt{\trace(A^TA)}$ be its trace norm.
%and $\|A\| = \max\{|Ax|: |x|=1\}$ be the operator norm.
%If $A$ is a systematic matrix ($A=A^T$), denote by $\l_{\min}(A)$ and $\l_{\max}(A)$
i%ts smallest and largest eigenvalue, respectively. By $A \le 0$ and $A<0$,
%we mean $A$ is non-negative and negative definite, respectively.
If both $a, b$ are real numbers, then $a\ve b = \min\{a, b\}$ and
$a\we b=\max\{a,b\}$.

Let $(\W ,{\cal F}, \{{\cal F}_t\}_{t\ge 0}, \PP)$ be a complete probability
space with a filtration  $\{{\cal F}_t\}_{t\ge 0}$ satisfying the
usual conditions (i.e. it is right continuous and ${\cal F}_0$
contains all $\PP$-null sets). If $A$ is a subset of $\W$, denote by $I_A$ its
indicator function; that is $I_A(\w)=1$ when $\w\in A$ and $0$ otherwise.
Let $B(t)  = (B_1(t),\cdots,B_m(t))^T$
be an $m$-dimensional
Brownian motion defined on the probability space.
Let $r(t)$, $t\ge 0$, be a right-continuous Markov chain on the same probability space
taking values in a finite state space $S=\{1, 2, \cdots, N\}$ with its generator
$\G=(\g_{ij})_{N\K N}$ given by
$$
\PP\{r(t+\D)=j | r(t)=i\}=
\begin{cases}
 \g_{ij}\D + o(\D) & \hbox{if }
i\not= j, \\
        1+\g_{ii}\D + o(\D) & \hbox{if }  i=j,
        \end{cases}
$$
where $\D>0$.  Here $\g_{ij} \ge 0$ is the transition rate from $i$ to $j$ if
$i\not=j$ while
$$
\g_{ii} = -\sum_{j\not= i} \g_{ij}.
$$
We assume that the Markov chain $r(\cdot)$ is independent of the Brownian motion
$B(\cdot)$.
Moreover,  denote by $M_{{\cal F}_t}(S)$ the
family of all ${\cal F}_t$-measurable $S$-valued
random variables.

Let $\T$ be a positive number.  Denote by $C([-\T,0]; R^n)$ the family of continuous functions
$\f:[-\T,0]\to R^n$ with the norm $\|\f\|=\sup_{-\T\le u\le 0}|\f(u)|$.
For $q>0$ and $t\ge 0$, denote by $L^p_{\F_t}(C)$ the family of
$\F_{t}$-measurable $C([-\T,0];R^n)$-valued random variables $\xi$ such that
$\E\|\xi\|^q <\8$, and by $L^q_{\F_t}(R^n)$ the family of
$\F_{t}$-measurable $R^n$-valued random variables $\eta$ such that
$\E|\eta|^q <\8$.

Consider the $n$-dimensional hybrid SDE (\ref{1.1}) on $t\ge 0$, where the coefficients
$f : R^n\K S\K R_+ \to R^n$ and
 $g: R^n\K S\K R_+ \to R^n \to R^{n\K n}$ are Borel measurable.
Assuming SDE (\ref{1.1}) is unstable, our aim here is to
  design a Borel measurable control function $u: R^n\K S\K R_+ \to R^n$
  so that the delay
feedback control $u(x(t-\tau)), r(t),t)$ will make the controlled
hybrid system (\ref{1.3}) become stable. Noting that (\ref{1.3}) is a
hybrid stochastic differential delay equation (SDDE), we naturally impose the
initial data
\begin{equation}\label{ID}
  \{x(\theta):-\T\le u\le 0\} = \f \in C([-\T,0];R^n)
  \quad \hbox{and} \quad
  r(0)=r_0\in S.
\end{equation}
This means at the current time $t=0$  the historical data of the state  $\{x(\theta):-\T\le u\le 0\}$ and the mode
$r(0)$ are available.  For the controlled SDDE (\ref{1.3}) to have a unique solution on $t\ge 0$ with the initial data (\ref{ID}), we impose the
global Lipschitz condition (see, e.g., \cite{M91,M94,M97,Moh}).

\begin{assp} \label{A2.1}
There there three positive constants $L_1, L_2$ and $L_3$ such that
\begin{eqnarray} \label{2.3a}
& & |f(x,i,t)-f(y,i,t)|\le L_1|x-y|, \nonumber  \\
& & |u(x,i,t)-u(y,i,t)|\le L_2|x-y|, \\
& & |g(x,i,t)-g(y,i,t)|\le L_3|x-y| \nonumber
 \end{eqnarray}
 for all $(x,y,i,t)\in R^n\K R^n\K S\K R_+$.  Moreover,
\begin{equation} \label{2.3b}
f(0,i,t)=0, \quad g(0,i,t)=0, \quad u(0,i,t) =0
\end{equation}
for all $(i,t)\in S\K R_+$.
\end{assp}

We see that this assumption implies the linear growth condition
\begin{equation} \label{2.3}
|f(x,i,t)| \le L_1|x|, \  |u(x,i,t)| \le L_2|x|, \   |g(x,i,t)| \le L_3|x|
\end{equation}
for all $(x,i,t)\in R^n\K S\K R_+$.
It is known (see, e.g., \cite{MY06}) that under Assumption \ref{A2.1},
the controlled SDDE (\ref{1.3})) with the initial data (\ref{ID}) has a unique solution $x(t)$ on $t\ge 0$ and the solution has the property that
\begin{equation} \label{qbd}
\E\|x_t\|^q < \8 \quad\hbox{for all } t \ge 0 \hbox{ and any } q>0,
\end{equation}
 where throughout this paper we use the notation $x_t=\{x(t+u): -\T\le u \le 0\}$ which is a $C([-\T,0];R^n)$-valued stochastic process on $t\ge 0$.   To emphasise the role of the initial data at time 0, we will denote
the solution by $x(t;\f,r_0,0)$ and the Markov chain by $r(t;r_0,0)$.
Let $t_0 \ge 0$, $x_{t_0}= \{x(t+u;\f,r_0,0): -\tau\le u\le 0\}$ and $r(t_0)=r(t_0;r_0,0)$.  Moreover, denote the unique solution of the SDDE (\ref{1.3}) on $t\ge t_0$ with the initial data $x_{t_0}$ and $r(t_0)$ at time $t_0$ by $x(t;x_{t_0},r(t_0),t_0)$ and the corresponding Markov chain by
$r(t;r(t_0),t_0)$.  We then see   the flow property that
 \begin{equation} \label{2.5}
 x(t;\f,r_0,0) = x(t;x_{t_0},r(t_0),t_0)
\ \hbox{and} \ r(t;r_0,0) = r(t;r(t_0),t_0)
\end{equation}
for all $t\ge t_0$.

Let us now return to the (non-delay) controlled hybrid SDE (\ref{1.2}).  Instead of $x(t)$ we will use $y(t)$ for the state to distinguish it from the solution of the SDDE (\ref{1.3}).  That is, we  consider the auxiliary controlled hybrid SDE
\begin{equation} \label{2.6}
dy(t)=  \big(f(y(t),r(t),t)+ u(y(t),r(t),t) \big) dt + g(y(t),r(t),t)dB(t).
\end{equation}
From now on, we will fix a number $p > 0$.
It is   known (see, e.g., \cite{MY06}) that under Assumption \ref{A2.1},
the SDE (\ref{2.6}) with the initial data $y(t_0)\in L^p_{{\cal F}_{t_0}}(R^n)$
and $r(t_0)=r_0\in M_{{\cal F}_{t_0}}(S)$ at time $t_0$ has a unique solution $y(t)$ on $t\ge t_0$ which has the property that $\E|y(t)|^p <\8$
for all $t\ge t_0$.   We will denote the solution by $y(t;y(t_0),r(t_0),t_0)$.
 As we mentioned before, there are already many papers devoted to the designation of the
control function $u: R^n\K S \K R_+ \to R^n$ for
this SDE to be exponentially stable in the $p$th moment.
We can therefore simply assume the exponential stability of this SDE.

\begin{assp} \label{A2.2}
Let $p>0$.  Assume that there is a pair of positive constants $M$ and $\g$ such that
the solution of the auxiliary
controlled hybrid SDE (\ref{2.6}) satisfies
\begin{equation} \label{2.7}
\E|y(t;y(t_0),r(t_0),t_0)|^p \le M \E|y(t_0)|^p e^{-\g(t-t_0)} \ \ \forall{t\ge t_0\ge 0}
\end{equation}
for all $y(t_0)\in L^p_{{\cal F}_{t_0}}(R^n)$
  and $r(t_0)\in M_{{\cal F}_{t_0}}(S)$.
\end{assp}

Our aim in this paper is to show that this same control function also makes
 the delay controlled system (\ref{1.3}) to be exponentially stable in the $p$th moment as long as
 $\T$ is sufficiently small (namely we make state observations
 frequently enough).  To prove this result, let us present a number of lemmas in the next section.

\section{Lemmas}

In this section, we will fix the initial data (\ref{ID}) arbitrarily.
We will write the solution $x(t;\f,r_0,0)=x(t)$ of the controlled hybrid SDDE (\ref{1.3}) with the initial data (\ref{ID}) and the Markov chain $r(t;r_0,0)=r(t)$ on $t\ge 0$. We emphasize that $x_t\in  L^p_{\F_t}(C)$
and $x(t)\in L^p_{\F_t}(R^n)$ (please recall (\ref{qbd})).
We also emphasize once again that we fix $p>0$ throughout this paper.

\begin{lemma}\label{L3.1}
Under Assumption, for any $t_0\ge 0$ and $T \ge 0$,
\begin{equation}\label{L3.1a}
   \sup_{t_0\le t\le t_0+T+\T} \E|x(t)|^p \le   K_1  \E\|x_{t_0}\|^p,
\end{equation}
\begin{equation}\label{L3.1b}
    \E\Big(\sup_{t_0\le t\le t_0+T+\T} |x(t)|^p\Big) \le   K_2  \E\|x_{t_0}\|^p,
\end{equation}
\begin{equation}\label{L3.1c}
   \sup_{t_0\le t\le t_0+T} \E\Big(\sup_{0\le u\le\tau} |x(t+u)-x(t)|^p\Big) \le   K_3  \E\|x_{t_0}\|^p,
\end{equation}
where
\begin{align}
K_1&= K_1(p,\T,T) = (1+\T)^{(1\we(0.5p))}
   e^{ p(T+\T)(L_1+L_2+0.5L_3^2[1\ve(p-1)] )},  \label{L3.1a1}  \\
K_2&= K_2(p,\T,T) \nonumber \\
& =
\left\{
\begin{array}{ll}
 4^{p-1} K_1(p,\T,T)  \big( 1 + (L_1^p+L_2^p)(T+\T)^p    + C_p L_3^p (T+\T)^{0.5p} \big), & \hbox{if } p\ge 2,  \\
  \big(  4 K_1(2,\T,T)\big[ (L_1^2+L_2^2) (T+\T)^2 +L_3^2(T+\T)\big] \big)^{p/2} , & \hbox{if } p\in (0,2),
 \end{array}
 \right.   \label{L3.1b1} \\
 K_3&= K_3(p,\T,T)=
\left\{
\begin{array}{ll}
 3^{p-1} K_1(p,\T,T)   [ (L_1^p+L_2^p)\T^p    + C_p L_3^p \T^{0.5p} ] , & \hbox{if } p\ge 2,  \\
 \big( 3 K_1(2,\T,T)  [(L_1^2+L_2^2)\T^2+4L_3^2\T]\big)^{p/2}, & \hbox{if } p\in (0,2),
 \end{array}
 \right.   \label{L3.1c1}
\end{align}
in which $C_p=[p^{p+1}/2(p-1)^{p-1}]^{p/2}$ for $p\ge 2$.  (Of course, $K_1$--$K_3$ depend on $L_1$--$L_3$ but we do not want to emphasise them explicitly as $L_1$--$L_3$ are fixed once the underlying SDE is given but we usually need to
choose $p, \T, T$ to fit into the underlying situation.)
\end{lemma}

\noindent
{\it Proof}.  We fix  $t_0\ge 0$ and $T \ge 0$ arbitrarily. We first prove the first assertion for $p\ge 2$. By the It\^o formula and  Assumption \ref{A2.1},  it is straightforward to show from (\ref{1.3}) that
for $t\in [t_0,t_0+T+\T]$
\begin{align*}
& \E|x(t)|^p \\
\le & \E|x(t_0)|^p +
  \E \int_{t_0}^t    \Big( p \big[ L_1|x(s)|^p+L_2|x(s)|^{p-1}|x(s-\T)|\big] + 0.5p(p-1)L_3^2|x(s)|^p \Big) ds   \\
\le & \E|x(t_0)|^p +   \big[pL_1+(p-1)L_2 + 0.5 p(p-1)L_3^2\big]  \int_{t_0}^t \E|x(s)|^p ds
+ L_2\int_{t_0}^t  \E |x(s-\T)|^p   ds.
\end{align*}
But
$$
\int_{t_0}^t  \E |x(s-\T)|^p   ds
\le \T \E \|x_{t_0}\|^p +  \int_{t_0}^t  \Big( \sup_{t_0\le u\le s} \E|x(u)|^p \Big) ds.
$$
Therefore
\begin{align*}
\E|x(t)|^p   \le  (1+\T)\E \|x_{t_0}\|^p +
   p\big[L_1+L_2 + 0.5 (p-1)L_3\big]  \int_{t_0}^t  \Big( \sup_{t_0\le u\le s} \E|x(u)|^p \Big) ds.
\end{align*}
As the last term on the right-hand-side of the inequality above is increasing in $t$,
we must have
$$
\sup_{t_0\le u\le t} \E|x(u)|^p \le  (1+\T)\E \|x_{t_0}\|^p +
   p\big[L_1+L_2 + 0.5 (p-1)L_3\big]  \int_{t_0}^t  \Big( \sup_{t_0\le u\le s} \E|x(u)|^p \Big) ds.
$$
An application of the well-known Gronwall inequality yields
$$
   \sup_{t_0\le u\le t_0+T} \E|x(u)|^p  \le  (1+\T) e^{p(T+\T)[L_1+L_2 + 0.5 (p-1)L_3] } \E\|x_{t_0}\|^p.
$$
This is the required assertion (\ref{L3.1a}) when $p\ge 2$.
When $p\in (0,2)$, we can apply the It\^o formula to $|x(t)|^2$ and
then take the conditional expectation given $\F_{t_0}$ to get that
\begin{align}\label{L3.1M}
 \E(|x(t)|^2|\F_{t_0})  \le  (1+\T) e^{2(T+\T)(L_1+L_2+0.5L_3^2)} \|x_{t_0}\|^2
\end{align}
for all $t\in [t_0,t_0+T+\T]$. Hence
\begin{align*}
 &  \E( |x(t)|^p |\F_{t_0}) \le \big(  \E( |x(t)|^2 |\F_{t_0}) \big)^{p/2}
 \le    (1+\T)^{p/2} e^{p(T+\T)(L_1+L_2+0.5L_3^2)} \|x_{t_0}\|^p.
\end{align*}
Taking the expectation on both sides gives the  required assertion (\ref{L3.1a}) when $p\in (0,2)$.  In other words, we have shown that (\ref{L3.1a}) holds for all $p>0$.

Let us proceed to prove the second assertion, namely (\ref{L3.1b}) for $p\ge 2$.
It is easy to show from (\ref{1.3}) that
\begin{align*}
& \E\Big(\sup_{t_0\le t\le t_0+T+\T} |x(t)|^p\Big) \\
  \le &  4^{p-1} \Big( \E|x(t_0)|^p +
(T+\T)^{p-1} \int_{t_0}^{t_0+T+\T} (L_1^p\E|x(s)|^p+L_2^p\E|x(s-\T)|^p) ds \\
  + & C_p  (T+\T)^{0.5p-1} L_3^p \int_{t_0}^{t_0+T+\T}  \E|x(s)|^p ds \Big),
\end{align*}
where $C_p$ has been defined in the statement of the lemma. Substituting (\ref{L3.1a})
into this yields (\ref{L3.1b}) when $p\ge 2$.  For $p\in (0,2)$, we note that
\begin{align*}
& \E\Big(\sup_{t_0\le t\le t_0+T+\T} |x(t)|^2\big|\F_{t_0}\Big) \\
  \le & 4 \Big( |x(t_0)|^2 +
(T+\T) \int_{t_0}^{t_0+T+\T} \big( L_1^2\E(|x(s)|^2|\F_{t_0})+ L_2^2 \E(|x(s-\T)|^2|\F_{t_0})\big) ds \\
  + & 4 L_3^2 \int_{t_0}^{t_0+T+\T}  \E(|x(s)|^p |\F_{t_0})  ds \Big),
\end{align*}
Substituting (\ref{L3.1M}) into the above implies
\begin{align*}
\E\Big(\sup_{t_0\le t\le t_0+T+\T} |x(t)|^2\big|\F_{t_0}\Big)
& \le 4 K_1 \big[ (L_1^2+L_2^2)(T+\T)^2 + L_3^2(T+\T) \big] \|x_{t_0}\|^2.
\end{align*}
Then
\begin{align*}
& \E\Big(\sup_{t_0\le t\le t_0+T+\T} |x(t)|^p\big|\F_{t_0}\Big)
\le  \Big[\E\Big(\sup_{t_0\le t\le t_0+T+\T} |x(t)|^2\big|\F_{t_0}\Big)\Big]^{p/2}
\\
& \le \big(  4 K_1 \big[ (L_1^2+L_2^2)(T+\T)^2 + L_3^2(T+\T) \big]  \big)^{p/2} \|x_{t_0}\|^p.
\end{align*}
Taking the expectation on both sides and recalling the definition of $K_2(p,\T,T)$, we see that the required assertion (\ref{L3.1b}) holds for $p\in (0,2)$ as well.

Similarly, we can show the third assertion (\ref{L3.1c}). The proof is complete. $\Box$

\begin{lemma} \label{L3.2}
Let Assumption \ref{A2.1} hold. Fix $t_0\ge \T$ and $T\ge 0$ arbitrarily.
Write $y(t;x(t_0),r(t_0),t_0)=y(t)$ for $t\ge t_0$.  Then, for $t\in [t_0,t_0+T+\T]$,
\begin{equation}\label{L3.2a}
  \E|x(t)-y(t)|^p \le  K_4(p,\T,L,T) \E\|x_{t_0}\|^p
\end{equation}
where
\begin{align}\label{L3.2b}
  K_4 &= K_4(p,\T,T) \nonumber  \\
 & =
\left\{
\begin{array}{ll}
 L_2 (T+\T) K_3(p,\T,T) e^{[ pL_1+(2p-1)L_2+0.5p(p-1)L_3^2] (T+\T) } ,
 & \hbox{if } p\ge 2,  \\
  \big( L_2 (T+\T) K_3(2,\T,T) e^{[ 2L_1+3L_2+L_3^2] (T+\T) }  \big)^{p/2} , & \hbox{if } p\in (0,2),
 \end{array}
 \right.
\end{align}
in which $K_3(p,\T,T)$ has been defined in Lemma \ref{L3.1}.
\end{lemma}

\noindent
{\it Proof}.  We first show the assertion for $p\ge 2$. By the It\^o formula and Assumption \ref{A2.1}, it is straightforward to show that for  $t\in [t_0,t_0+T+\T]$,
\begin{align*}
  \E|x(t)-y(t)|^p \le & \E \int_{t_0}^t \Big(  (pL_1+0.5p(p-1)L_3^2)|x(s)-y(s)|^p
  \\
 &  + pL_1|x(s)-y(s)|^{p-1} |x(s-\T)-y(s)| \Big) ds
   \end{align*}
But
\begin{align*}
 & p|x(s)-y(s)|^{p-1} |x(s-\T)-y(s)| \\
  \le &
  p|x(s)-y(s)|^p + p|x(s)-y(s)|^{p-1} |x(s-\T)-x(s)| \\
 \le & (2p-1)|x(s)-y(s)|^p + |x(s)-x(s-\T))|^p
\end{align*}
Hence
\begin{align*}
  \E|x(t)-y(t)|^p \le & \big[ pL_1+(2p-1)L_2+0.5p(p-1)L_3^2\big] \int_{t_0}^t \E |x(s)-y(s)|^p ds
  \\
+  &    L_2\int_{t_0}^t \E|x(s-\T)-x(s)|^p ds  \\
\le & \big[ pL_1+(2p-1)L_2+0.5p(p-1)L_3^2\big]  \int_{t_0}^t \E |x(s)-y(s)|^p ds \\
  + & L_2 (T+\T) K_3 \E\|x_{t_0}\|^p.
   \end{align*}
An application of the Gronwall inequality gives the assertion for $p\ge 2$.
Let us now consider the case when $p\in (0,2)$.  In a similar way as Lemma \ref{L3.1} was proved, we can show that
\begin{align*}
  \E(|x(t)-y(t)|^2 |\F_{t_0})
  \le  L_2 (T+\T) K_3(2,\T,T) e^{[ 2L_1+3L_2+L_3^2] (T+\T) }   \|x_{t_0}\|^2
   \end{align*}
for $t\in [t_0,t_0+T+\T]$.  We can then show that the assertion holds for $p\in (0,2)$ using the technique of conditional expectation as we did in the proof of Lemma \ref{L3.1}.  The proof is therefore complete.

\section{Main Result}

We can now form our main theorem in this paper.

\begin{theorem}\label{T4.1}
Under Assumptions \ref{A2.1} and \ref{A2.2}, there is a positive number $\T^*$ such that the solution of the controlled SDDE (\ref{1.3}) has the properties that
\begin{equation}\label{T4.1a}
\limsup_{t\to\8}  \frac{1}{t} \log(\E|x(t;\f,r_0,0)|^p) < 0
\end{equation}
and
\begin{equation}\label{T4.1b}
\limsup_{t\to\8}  \frac{1}{t} \log( |x(t;\f,r_0,0)| ) < 0 \quad a.s
\end{equation}
 as long as $\T < \T^*$.  In other words, the controlled SDDE (\ref{1.3})
is exponentially stable in the $p$th moment as well as with probability one
provided  $\T < \T^*$.
\end{theorem}

Before the proof, let us make an important remark on how to determine the value of $\T^*$ so that this theorem can be implemented in practice.

\begin{rmk} \label{R4.2}
The use of this theorem in practice depends very much on the value of $\T^*$.
We describe a method to determine it.
Set $p_0=0\ve(p-1)$. Choose a constant $\e\in (0,1)$ and let
\begin{equation} \label{T4.1c}
T= \frac{1}{\g} \log\Big( \frac{2^{2p_0} M}{\e}\Big).
\end{equation}
Let $\T^*$  be the unique positive root to the following equation
\begin{equation} \label{startau}
  2^{p_0}  [ 2^{p_0}K_4(p,\T,T) +  K_3(p,\T,T)  ] = 1-\e
\end{equation}
of $\T$, where both  $K_3(p,\T,T)$ and $K_4(p,\T,T)$ have been defined in Section 3.  We observe that $\T^*$ exists uniquely and is positive as the left-hand-side term of equation (\ref{startau}) is an increasing continuous function of $\T$ which starts from 0 when $\T=0$ and tends to infinity as $\T\to\8$.   However, we do not have the explicit formula for the root $\T^*$ though it can be solved numerically, for example, by MATHEMATICA.
We also observe that it is more desirable in practice if we could find a larger value of $\T^*$.
Note that once $p, L,M,\g$ are given, the root $\T^*$  depends on the choice of $\e$.
That is, $\T^*=\T^*(\e)$.  It would be useful if we could find the optimal $\bar \e \in (0,1)$ in the sense that
$$
\T^*(\bar\e) = \sup_{\e\in (0,1)} \T^*(\e).
$$
However, this is an open problem.
\end{rmk}

\noindent
{\it Proof of Theorem \ref{T4.1}}.  To make it clearer, we divide the proof
into three steps.

\emph{Step 1}. We will simply write $K_3(p,\T,T) = K_3$ and $K_4(p,\T,T)=K_4$.  We let $\T^*$ be determined in the way as described in Remark \ref{R4.2}.
Fix $\T \in (0, \T^*)$ and the initial data (\ref{ID}).
Write $x(t;\f,r_0,0)=x(t)$ and $r(t;r_0,0)=r(t)$ for $t\ge 0$.
Let us first consider $x(t)$ on $t\in [\T, 2\T+T]$ which can be regarded as the solution of the  SDDE (\ref{1.3}) with initial data
$x_\T$ and $r(\T)$ at time $t=\T$.
Also consider the solution  $y(t; x(\T),r(\T),\T)$ of
the SDE (\ref{2.6}) on $t\in [\T, \T+T]$ with initial data
$x(\T)$ and $r(\T)$ at time $t=\T$. We simply write
$y(\T+T; x(\T),r(\T),\T)=y(\T+T)$.    By Assumption \ref{A2.2},
\begin{equation}\label{T4.1d}
\E|y(\T+T)|^p \le M e^{-\g T}  \E|x(\T)|^p.
\end{equation}
Moreover, by the elementary inequality $(a+b)^p\le 2^{p_0}(a^p+b^p)$ for any $a,b\ge 0$
(please recall that $p_0:=0\ve (p-1)$ which has been defined in Remark \ref{R4.2}),
we have
$$
\E|x(\T+T)|^p \le 2^{p_0} \Big( \E|y(\T+T)|^p +
 \E|x(\T+T)-y(\T+T)|^p \Big).
$$
Using (\ref{T4.1d})
and Lemma \ref{L3.2},
we get
\begin{equation} \label{T4.1e}
\E|x(\T+T)|^p \le 2^{p_0}
\Big(M e^{-\g T}  \E|x(\T)|^p +  K_4 \E\|x_\T|^p \Big)
\le 2^{p_0}(M e^{-\g T} +K_4) \E\|x_\T\|^p.
\end{equation}
On the other hand, by Lemma \ref{L3.1}, we have
\begin{align} \label{T4.1f}
  \E\|x_{2\T+T}\|^p
\le & 2^{p_0} \E|x(\T+T)|^p
+2^{p_0} \E\Big( \sup_{0\le u\le \T} |x(\T+T)-x(\T+T+u)|^p \Big)
\nonumber \\
\le & 2^{p_0} \E|x(\T+T)|^p
+2^{p_0} K_3 \E\|x_\T\|^p.
\end{align}
Substituting (\ref{T4.1e}) into (\ref{T4.1f}) and noting from (\ref{T4.1c})
that
$$
 2^{2p_0} M e^{-\g T} = \e,
$$
we get
\begin{align} \label{T4.1g}
  \E\|x_{2\T+T}\|^p
\le  \big[ \e + 2^{p_0}  ( 2^{p_0}K_4+  K_3  ) \big] \E\|x_\T\|^p.
\end{align}
But, as $\T <\T^*$, we see from (\ref{startau}) that
\begin{align*}
\e + 2^{p_0}  ( 2^{p_0}K_4+  K_3  ) < 1.
\end{align*}
We may therefore write
$$
\e + 2^{p_0}  ( 2^{p_0}K_4+  K_3  ) = e^{-\l (\T+T)}
$$
for some $\l > 0$. It then follows from (\ref{T4.1g}) that
\begin{equation} \label{T4.1h}
\E\|x_{2\T+T}\|^p \le e^{-\l (\T+T)} \E\|x_\T\|^p.
\end{equation}

\emph{Step 2}.
Let us now consider the solution $x(t)$ on $t\in [2\T +T, \T+2(\T+T)]$. By property (\ref{2.5}), this can be regarded as the solution of the SDDE (\ref{1.3}) with the initial data $x_{2\T+T}$ and $r(2\T+T)$ at $t=2\T+T$. In the same way as we did in Step 1, we can show
$$
\E\|x_{\T+2(\T+T)}\|^p \le e^{-\l (\T+T)} \E\|x_{2\T+T}\|^p.
$$
This, together with (\ref{T4.1i}),  implies
$$
\E\|x_{\T+2(\T+T)}\|^p \le e^{- 2\l (\T+T)} \E\|x_\T\|^p.
$$
Repeating this procedure, we have
\begin{equation} \label{T4.1i}
\E\|x_{\T+k(\T+T)}\|^p \le e^{- k \l (\T+T)} \E\|x_\T\|^p
\end{equation}
for all $k=1,2,\cdots$.  But, by Lemma \ref{L3.1}
$$
\E \Big( \sup_{\T+k(\T+T)\le t\le \T+(k+1)(\T+T)} |x(t)|^p \Big) \le
K_2 \E\|x_{\T+k(\T+T)}\|^p.
$$
This, together with (\ref{T4.1i}), yields
\begin{equation}\label{T4.1j}
\E \Big( \sup_{\T+k(\T+T)\le t\le \T+(k+1)(\T+T)} |x(t)|^p \Big) \le
K_2 e^{- k \l (\T+T)} \E\|x_\T\|^p
\end{equation}
for all $k=0,1,2,\cdots$.  Hence, for $t\in [\T+k(\T+T), \T+(k+1)(\T+T)]$,
$$
\frac{1}{t}\log(\E|x(t)|^p) \le
\frac{\log(K_2\E\|x_\T\|^p) -   k \l (\T+T)}{\T+k(\T+T)}.
$$
This implies
$$
\limsup_{t\to\8} \frac{1}{t}\log(\E|x(t)|^p) \le -\l.
$$
In other words, we have shown the required assertion (\ref{T4.1a}).

\emph{Step 3}. It now follows from (\ref{T4.1j}) that
$$
\PP\Big( \sup_{\T+k(\T+T)\le t\le \T+(k+1)(\T+T)} |x(t)|^p
\ge e^{- 0.5 k \l (\T+T)} \Big)
\le K_2 e^{-0.5 k \l (\T+T)} \E\|x_\T\|^p
$$
for all $k\ge 1$.
By the Borel--Cantelli lemma (see, e.g., \cite[Lemma 2.4 on page 7]{M97}), we obtain that
for almost all $\w\in \W$, there is an
integer $k_0=k_0(\w)$ such that
$$
\sup_{\T+k(\T+T)\le t\le \T+(k+1)(\T+T)} |x(t)|^p
< e^{- 0.5 k \l (\T+T)} \quad \forall k\ge k_0(\w).
$$
This implies easily that
$$
\limsup_{t\to\8} \frac{1}{t} \log(|x(t,\w)|) \le -\frac{ \l } {2p}
$$
for almost all $\w\in\W$. The other assertion (\ref{T4.1b}) must therefore hold. The proof is hence complete.

\section{Implementation}

In this section we will demonstrate how to implement our theory in order to stabilise the given unstable hybrid SDE (\ref{1.1}) by a delay feedback control $u(x(t-\T),r(t),t)$ in the drift.  Our new Theorem \ref{T4.1} enables us to transfer the stabilisation problem of (\ref{1.3}) to the classical stabilisation problem of (\ref{2.6}) where the feedback control $u(x(t),r(t),t)$ is of no-delay.  The use of our new Theorem \ref{T4.1} depends on the design of the control function $u(x,i,t)$ that makes the
controlled SDE (\ref{2.6}) become exponentially stable in the $p$th moment
as described in Assumption \ref{A2.2}.  There are lots of known criteria on the $p$th moment exponential stability of hybrid SDEs (see, e.g., \cite{M99,MYY,YM04,YH}), which can be applied to design the control function.
What we are going to demonstrate here is to apply \cite[Theorem 5.8 on page 166]{MY06} to establish a criterion for the control function to satisfy.
For this purpose, we will impose a new assumption.

\begin{assp} \label{A5.1}
Let $p>0$. Assume that there are real numbers $\a_i$, $i\in S$, such that
\begin{align}\label{A5.1a}
  \frac{p}{|x|^2} \big(   x^T[f(x,i,t)+u(x,i,t)] + 0.5|g(x,i,t)|^2 \big) -\frac{ p(2-p)}{2|x|^4}|x^Tg(x,i,t)|^2 \le - \a_i
\end{align}
for all $(x,i,t)\in (R^n-\{0\})\K S\K R_+$ while
\begin{align}\label{A5.1b}
{\cal A}:=   \diag(\a_1, \cdots, \a_N) - \G
\end{align}
is a non-singular M-matrix.
\end{assp}

The following theorem shows that if the control function $u(x,t,i)$ makes Assumption \ref{A5.1} to hold, then the controlled SDE (\ref{2.6}) is exponentially stable in the $p$th moment.

\begin{theorem}\label{T5.2}
  Let Assumptions \ref{A2.1} and \ref{A5.1} hold.
  Then Assumption \ref{A2.2} holds with
  \begin{equation} \label{T5.2a}
  M = \beta_2/\beta_1 \quad\hbox{and}\quad    \g = 1/\beta_2,
  \end{equation}
  where
  \begin{equation} \label{T5.2b}
  (\theta_1,\cdots, \theta_N)^T = {\cal A}^{-1}(1,\cdots, 1)^T,
  \ \  \beta_1=\min_{i\in S}\theta_i, \ \  \beta_2=\max_{i\in S}\theta_i.
  \end{equation}
Consequently, Theorem \ref{T4.1} holds under Assumptions \ref{A2.1} and \ref{A5.1}.
\end{theorem}

\noindent
{\it Proof}.  We observe that all $\theta_i$'s are positive as ${\cal A}$ is a nonsingular M-matrix (see, e.g., \cite{MY06}).  It follows from (\ref{T5.2b}) that
\begin{equation} \label{T5.2c}
   \a_i\theta_i - \sum_{j=1}^N \g_{ij} \theta_j  = 1, \quad i\in S.
  \end{equation}

We will apply \cite[Theorem 5.8 on page 166]{MY06} to prove this theorem. We first consider the controlled SDE (\ref{2.6}) on $t\ge t_0$ in the case where the initial data are deterministic,
namely $y(t_0)=y_0\in R^n$ and $r(t_0)=r_0\in S$ at time $t_0 (\ge 0)$. We will write the solution $y(t;y_0,r_0,t_0)=x(t)$ and the Markov chain $r(t;r_0,t_0)=r(t)$.   Clearly, the assertion holds if $y_0=0$ so we need to consider $y_0\not= 0$.  In this case, $x(t)\not= 0$ a.s. for all $t\ge 0$ (see, e.g., \cite[Lemma 5.1 on page 164]{MY06}). Define the Lyapunov function $V: R^n-\{0\})\K S\K R_+$ by
$V(x,i,t)=\theta_i|x|^p$. So
$$
\beta_1|x|^p \le V(x,i,t) \le \beta_2 |x|^p, \quad (x,i,t) \in R^n-\{0\})\K S\K R_+,
$$
where both $\beta_1$ and $\beta_2$ have been defined in the statement of the theorem.
 Moreover,
the generalised It\^o formula (see, e.g., \cite[Theorem 1.45 on page 48]{MY06}) shows
$$
dV(x(t),r(t),t) = LV(x(t),r(t),t) dt + dM(t),
$$
where $M(t)$ is a local Martingale on $t\ge 0$ (but its explicit form is of no use here) and
\begin{align*}
LV(x,i,t) & = \theta_i |x|^p\Big( \frac{p}{|x|^2} \big(   x^T[f(x,i,t)+u(x,i,t)] + 0.5|g(x,i,t)|^2 \big) \\
& -\frac{ p(2-p)}{2|x|^4}|x^Tg(x,t,i)|^2 \Big)
     + \sum_{j=1}^N \g_{ij} \theta_j |x|^p
\end{align*}
for $(x,i,t) \in R^n-\{0\})\K S\K R_+$.
By (\ref{A5.1a}) and (\ref{T5.2c}), we then have
$$
LV(x,i,t)  \le -\Big( \a_i\theta_i - \sum_{j=1}^N \g_{ij} \theta_j \Big) |x|^p = -|x|^p.
$$
An application of \cite[Theorem 5.8 on page 166]{MY06} yields
$$
\E|x(t)|^p \le M |y_0|^p e^{-\g (t-t_0)}, \quad t\ge t_0,
$$
where both $M$ and $\g$ have been defined in the statement of Theorem \ref{T5.2}.
We now consider the general case, namely the controlled SDE (\ref{2.6}) on $t\ge t_0$ with the initial data
 $y(t_0)=y_0\in L^p_{\F_{t_0}}(R^n)$ and $r(t_0)=r_0\in M_{\F_{t_0}}(S)$ at time $t_0$. In this case, by the technique of conditional expectation, we derive
 $$
 \E|x(t)|^p = \E\big( \E(|x(t)|^p |\F_{t_0} )\big)
\le \E\big( M |y_0|^p e^{-\l (t-t_0)} \big) = M \E|y_0|^p e^{-\g (t-t_0)}
 $$
for $t\ge t_0$.  In other words, we have shown that Assumption \ref{A2.2} holds.  Consequently, Theorem \ref{T4.1} holds under Assumptions \ref{A2.1} and \ref{A5.1}.  The proof is complete.

Accordingly, we can implement our theory in two steps assuming that the coefficients $f$ and $g$ of the given hybrid SDE (\ref{1.1}) satisfy Assumption \ref{A2.1}:

\begin{itemize}
\item[Step 1] Design the control function $u(x,i,t)$ which satisfies Assumptions \ref{A2.1} and \ref{A5.1}.  Compute $\theta_i$'s by (\ref{T5.2b}) and determine both $M$ and $\g$ by (\ref{T5.2a}).

\item[Step 2] Choose a constant $\e\in (0,1)$ and compute $T$ by (\ref{T4.1c}). Find the unique positive root $\T^*$ of equation (\ref{startau}) numerically. Make sure the delay feedback control $u(x(t-\T),r(t),t)$ acts quickly enough in the sense $\T < \T^*$. Then the controlled hybrid SDDE (\ref{1.3}) is exponentially stable in the $p$th moment as well as in probability one.
\end{itemize}

Let us discuss an example to illustrate our theory.

\begin{expl} \label{E5.3}
{\rm
Consider a second order SDE
$$
\ddot{z}(t)+(a_{r(t)}+b_{r(t)}\dot{B}(t))\dot{z}(t)
+z(t)+c_{r(t)}\sin(z(t)) = 0,
$$
where $\dot{B}(t)$ is a scalar white noise
(informally thought as the derivative of a scalar Brownian motion $B(t)$), $r(t)$ is a Markov chain
taking values in the state space $S=\{1, 2\}$ with the generator
$$
\G = \begin{pmatrix}
       -1& 1    \\
       2 & -2
     \end{pmatrix},
$$
and the coefficients are specified by
$$
a_1=0.5, \ a_2=0.1, \ b_1=0.4, \ b_2=0.5, \
c_1=0.1, \ c_2=-0.1.
$$
This SDE has been used to describe, for example, the nonlinear hybrid stochastic oscillator (see, e.g., \cite{M97}).
Introducing $x(t)=(x_1(t),x_2(t))^T = (z(t),\dot{z}(t))^T$, we can write the oscillator
as the two-dimensional hybrid SDE
\begin{equation}\label{E5.3a}
  dx(t) = f(x(t),r(t)) dt + g(x(t),r(t)) dB(t),
\end{equation}
where
$$
f(x,i) = \begin{pmatrix}
        x_2   \\
        -x_1-c_i\sin(x_1)-a_ix_2
     \end{pmatrix}
 \quad\hbox{and}\quad
   g(x,i) = \begin{pmatrix}
        0  \\
       - b_i x_2
     \end{pmatrix}.
$$
The computer simulation (see Figure 5.1) shows this given hybrid SDE is unstable.
\begin{figure}[h!]
\begin{center}
\includegraphics[angle=0,height=10cm,width=12cm]{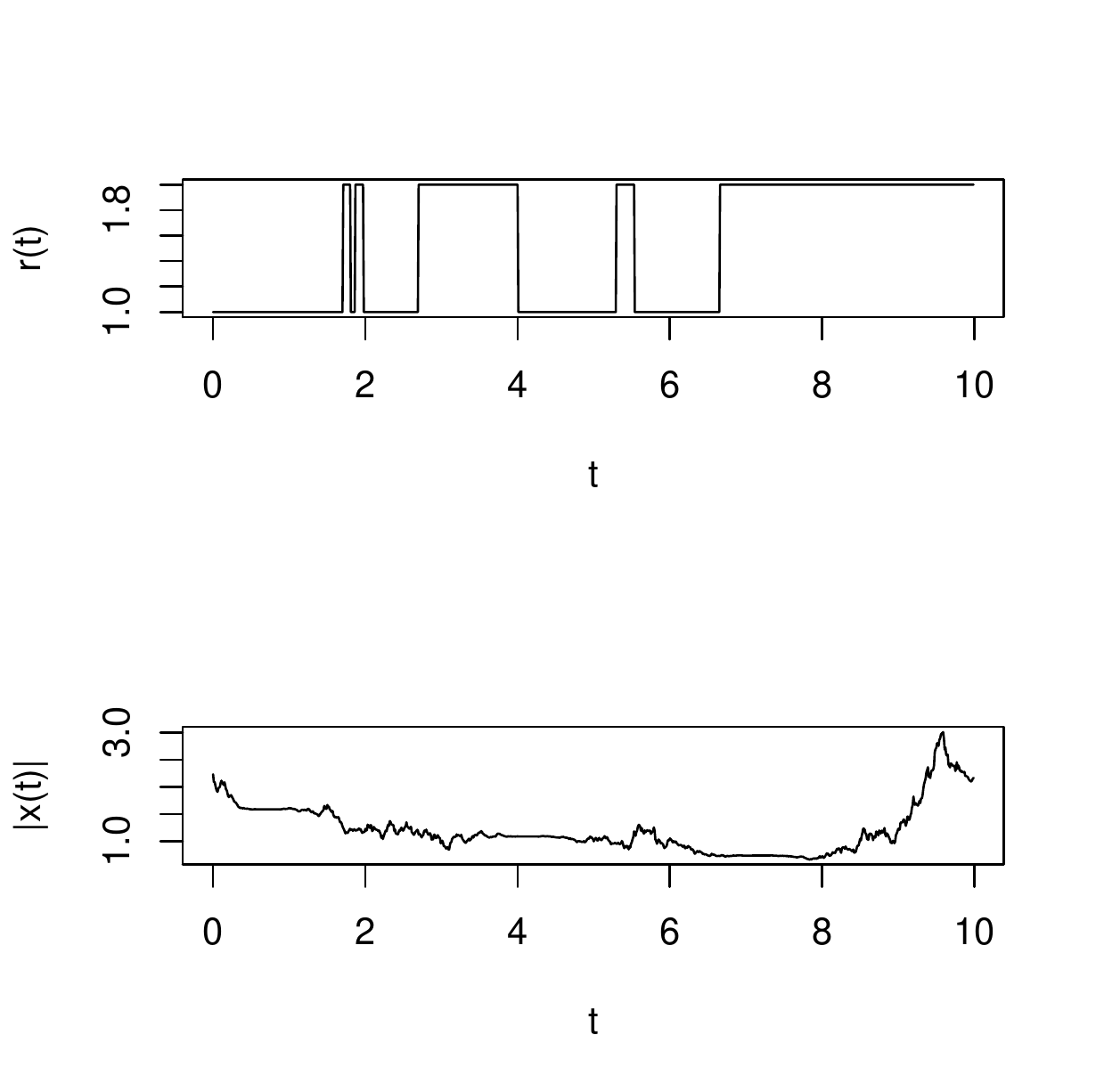}
 \end{center}
\begin{center}
\small{Figure 5.1: The computer simulation of
the sample paths of the Markov chain and the SDE  (\ref{E5.3a})
the initial data
$x_1(0)=1$, $x_2(0)=2$ and
$r(0)=1$
using the Euler--Maruyama method (see, e.g., \cite{LMY19}) with step size $10^{-5}$.}
\end{center}
\end{figure}

Let us now apply our new theory to design a delay feedback control to stabilise the SDE.  To show our theory can be applied to cope with various practical situations, we consider a structure feedback control in this example. Due to the page limit, we only discuss an interesting situation, where
\begin{itemize}
  \item the state, in both modes, could only be observed in $x_1$-component and the control could only be fed  into $x_1$-component too.
\end{itemize}
\noindent
For example, this is the case when $x_1$ represents the velocity
and $x_2$ the distance in a hybrid stochastic oscillation system while only the velocity is observable and controllable.

To make it simple, we will only seek for a linear control function. In terms of mathematics, our control function has the form
\begin{equation}\label{E5.3b}
u(x,i) =  \begin{pmatrix}
        -d_i x_1  \\
       0
     \end{pmatrix} \quad\hbox{for } (x,i)\in R^2\K S,
\end{equation}
where $d_1$ and $d_2$ are both positive numbers to be chosen.
It is straightforward to show that Assumption \ref{A2.1} is satisfied with
$$
L_1= 1.118034, \  L_2=d_1\ve d_2, \ L_3=0.5.
$$
It is also easy to show
\begin{align*}
 &  \frac{p}{|x|^2} \big(   x^T[f(x,i,t)+u(x,i,t)] + 0.5|g(x,i,t)|^2 \big) -\frac{ p(2-p)}{2|x|^4}|x^Tg(x,i,t)|^2 \le \nonumber \\
 & \quad \le \frac{p}{|x|^4} (x_1^2, x_2^2) Q_i (x_1^2,x_2^2)^T,
\end{align*}
for $(x,i)\in R^2-\{0\})\K S$, where
$$
Q_i =  \begin{pmatrix}
        |c_i|-d_i & 0.5(-a_i+0.25b_i^2-d_i) \\
       0.5(-a_i+0.25b_i^2-d_i) & |c_i|-a_i-0.5(1-p)b_i^2
     \end{pmatrix},
$$
namely
$$
Q_1 =  \begin{pmatrix}
        0.1-d_1 & -0.23-0.5d_1 \\
       -0.23-0.5d_1 & -0.464+0.08p
     \end{pmatrix}
     \ \ \hbox{and} \ \
 Q_2 =  \begin{pmatrix}
        0.1-d_2 & 0.28125-0.5d_2 \\
       0.28125-0.5d_2 & -0.125(1-p)
     \end{pmatrix}.
$$
In this example, we aim to stabilise the SDE in the sense of almost sure exponential stability so it is wise to choose  $p\in (0,1)$ to make use of the stabilisation effect of the Brownian motion
(see, e.g., \cite{MYY}).   We choose $p=0.99$ and $d_1$ for $0.1-d_1=-0.464+0.08p$, namely $d_1=0.4848$, while
$d_2$ for $0.28125-0.5d_2 = -0.125(1-p)$, namely $d_2=0.5650$. Consequently
$$
Q_1 = \begin{pmatrix}
        -0.3848 & -0.4724 \\
       -0.4724 &  -0.3848
     \end{pmatrix}
     \quad\hbox{and}\quad
 Q_2 =  \begin{pmatrix}
        -0.4650 & -0.0012 \\
       -0.0012 & -0.0013
     \end{pmatrix}.
$$
It is then easy to see that $\a_1$ and $\a_2$ in Assumption \ref{A5.1} are: $\a_1=0.3848$ and $\a_2=0.0012$.  The matrix defined by (\ref{A5.1b}) becomes
$$
{\cal A} = \begin{pmatrix}
        1.3848 & -1 \\
       -2 &  2.0012
     \end{pmatrix},
$$
which is a nonsingular M-matrix. In other words, we have verified Assumption \ref{A5.1} for $p=0.99$.  We can then determine $\theta_1=3.891286$, $\theta_2=4.388653$ and, hence,
$M= 1.127816$ and $\g=0.2278604$ by (\ref{T5.2b}) and (\ref{T5.2a}), respectively.  By Theorem \ref{T5.2}, the controlled SDE
\begin{equation}\label{E5.3c}
  dx(t) = [f(x(t),r(t))+u(x(t),r(t))] dt + g(x(t),r(t)) dB(t)
\end{equation}
with the control function $u(x,i)$ defined by (\ref{E5.3b})
is almost surely exponentially stable (please note that the $p$th moment exponential stability implies the almost sure  exponential  stability \cite{M97}).
The computer simulation (see Figure 5.2) supports this theoretical result clearly.
\begin{figure}[h!]
\begin{center}
\includegraphics[angle=0,height=10cm,width=12cm]{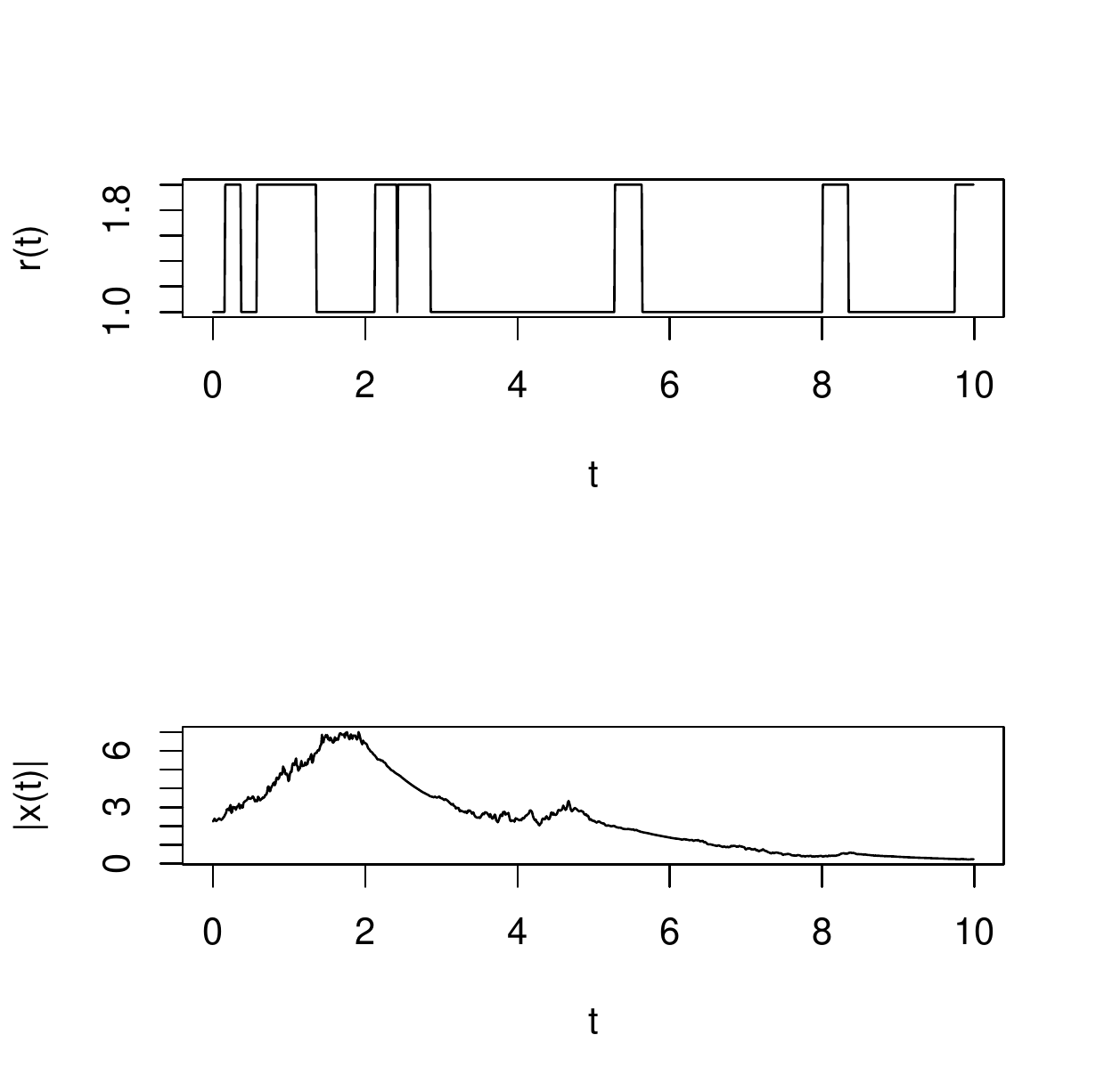}
 \end{center}
\begin{center}
\small{Figure 5.2: The computer simulation of
the sample paths of the Markov chain and the controlled SDE  (\ref{E5.3c}) the initial data
$x_1(0)=1$, $x_2(0)=2$ and
$r(0)=1$
using the Euler--Maruyama method with step size $10^{-5}$.}
\end{center}
\end{figure}

However, our aim is to use the delay feedback control.
For this purpose,
we further choose $\e=0.94$ and compute $T=0.7994283$ by (\ref{T4.1c}).
 Equation (\ref{startau}) becomes
\begin{equation}
K_4(0.99,\tau,0.7994283)+K_3(0.99,\tau,0.7994283)
=0.06,
\end{equation}
which has the unique positive root $\tau^*=2.93\K 10^{-6}$
(that is about 92 seconds if the time unit is of year).
%(a second ($=3.17*10^{-8}$ year)
By Theorem \ref{T5.2}, the delay-feedback controlled system
\begin{equation}\label{E5.3d}
  dx(t) = [f(x(t),r(t))+u(x(t-\tau),r(t))] dt + g(x(t),r(t)) dB(t)
\end{equation}
with the control function $u(x,i)$ defined by (\ref{E5.3b})
is almost surely exponentially stable
as long as $\tau < 2.93\K 10^{-6}$.
Once again, the computer simulation (see Figure 5.3) supports this theoretical result clearly.
\begin{figure}[h!]
\begin{center}
\includegraphics[angle=0,height=10cm,width=12cm]{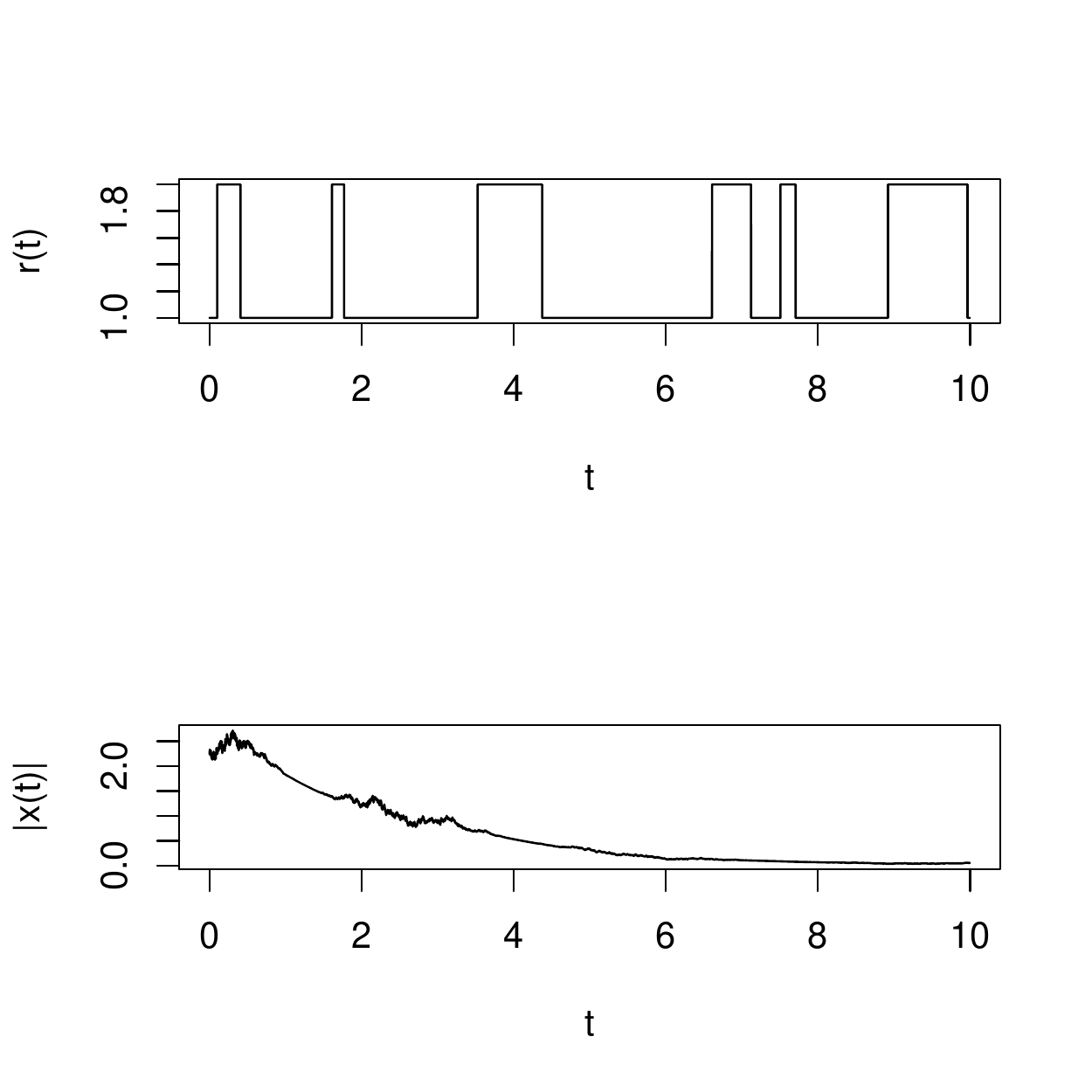}
 \end{center}
\begin{center}
\small{Figure 5.3: The computer simulation of
the sample paths of the Markov chain and the controlled system  (\ref{E5.3d}) with $\tau=10^{-6}$ and the initial data
$x_1(u)=1$, $x_2(u)=2$ for $u\in [-10^{-6},0]$ and
$r(0)=1$
using the Euler--Maruyama method with step size $10^{-7}$.}
\end{center}
\end{figure}
}\end{expl}

\begin{rmk}
Actually, the $\tau^*$ obtained in this example is optimal, to a certain degree according to our theory. More precisely, given that all the coefficients and $\Gamma$ in \eqref{E5.3a} are known, the largest $\tau^*$ is found by a numerical programme which searches for $p\in (0,1)$ and $\epsilon\in (0,1)$ to maximise $\tau^*$  according to Theorem \ref{T5.2}.
\end{rmk}

\section{Conclusion}

For some unstable hybrid stochastic differential equations, it is much harder to design the feedback control based on past states than current states. However, the feedback control based on past states are more practical than that based on current states. In this paper, we proposed a new theorem to connect the delay feedback control to the feedback control without delays. Such a result makes it possible to construct the delay feedback control $u(x(t-\tau),r(t),t)$, given that the feedback control $u(x(t),r(t),t)$ is known. Therefore, various existing results on the classical stabilisation problem together with the novel theorem proved in this paper enable us to design  the delay feedback control much more easily.
\par
Numerical simulations were provided to demonstrate the theoretical results as well as the way to find the lower bound of the length of the time delay. In addition, the optimal lower bound (according to our theory) was found numerically for a given unstable hybrid stochastic differential equation and a given format of the control function.

\section*{Appendix}

Consider a scalar SDE
\begin{equation} \label{7.1}
dx(t) = -x(t) dt + x^2(t)dB(t)
\end{equation}
on $t\ge 0$, where $B(t)$ is a scalar Brownian.   It is known (see, e.g., \cite{M97}) that for any given initial value $x(0)=x_0\in R$,
there is a unique global solution to the SDE but the equation is not exponentially stable in the $p$th moment for any $p\ge 2$.

We design a (non-delay) feedback control $-2x^3(t)$ in the drift so that the controlled SDE is
\begin{equation} \label{7.2}
dx(t) = (-x(t) -2 x^3(t)) dt + x^2(t)dB(t).
\end{equation}
It is almost straightforward to show that for any initial value $x(0)=x_0\in R$, the solution of this controlled SDE satisfies
\begin{equation} \label{7.3}
\E|x(t)|^4 \le |x_0|^4 e^{-4t}, \  \  \forall t\ge 0.
\end{equation}
 In other words, the feedback control stabilises the given system (\ref{7.1}) in the sense of the 4th moment exponential stability.

 Let us now show that for any $\e >0$ (no matter how small it is), the corresponding delay feedback control $-2x^3(t-\e)$ in the drift can NOT stabilise the given system (\ref{7.1}) in the sense of the 4th moment exponential stability.  In other words, we will show
that  the corresponding delay feedback controlled system
\begin{equation} \label{7.4}
dx(t) = (-x(t) - 2 x^3(t-\e)) dt + x^2(t)dB(t)
\end{equation}
is NOT exponentially stable in the 4th moment.  It is easy to see that this SDDE has a unique global solution $x(t)$
for any initial data $\{x(\o):-\e \le \o\le 0\}
=\f \in C([-\e, 0];R)$.  If this SDDE were exponentially stable in the 4th moment, then for any initial data $
=\f \in C([-\e, 0];R)$, the 4th moment $\E|x(t)|^4$ of the solution must be finite and tend to 0 exponentially fast as $t\to\8$.  We consider
a special initial data $\f \in C([-\e, 0];R)$ such that
\begin{equation} \label{7.5}
\f(0) = \bar z  \quad\hbox{and}\quad  8|\f(\o)|^3 \le \bar z^2 \ \hbox{ for } \o\in [-\e, -0.5\e],
\end{equation}
where $\bar z$ is the unique positive root to the following equation on $z\ge 2$
\begin{equation} \label{7.6}
0.5 - \frac{2}{z^2} = e^{-\e(2+0.5z^2)}.
\end{equation}
As we indicated that $\E|x(t)|^4<\8$ for all $t\ge 0$,
we can show easily that, for $t\in [0,0.5\e]$,
\begin{align} \label{7.7}
\E|x(t)|^2 &= |x(0)|^2 +\E\int_0^t \big( 2x(s)[-x(s)-2x^3(s-\e)] +|x(s)|^4 \big) ds \nonumber\\
& \ge  |x(0)|^2 +\int_0^t \big( - 2\E|x(s)|^2  -4 |\f(s-\e)|^3 \E|x(s)|  +\E|x(s)|^4 \big) ds
\nonumber\\
& \ge  \bar z^2 +\int_0^t \big( - 2\E|x(s)|^2  -0.5 \bar z^2 \E|x(s)|  +\E|x(s)|^4 \big) ds.
\end{align}
It is easy to see that  $\E|x(t)|^2 \ge \bar z^2 >4$ for $t\in [0,0.5\e]$.  By the H\"older inequality,
$$
\E|x(s)|^4 \ge (\E|x(s)|^2)^2 \quad\hbox{and}\quad
\E|x(s)| \le \sqrt{\E|x(s)|^2} \le \E|x(s)|^2
$$
for $s\in [0,0.5\e]$.     It then follows from (\ref{7.7}) that
\begin{align} \label{7.8}
\E|x(t)|^2
 \ge  \bar z^2 +\int_0^t \big[ - (2+0.5 \bar z^2) \E|x(s)|^2  +(\E|x(s)|^2)^2 \big] ds
\end{align}
for $t\in [0,0.5\e]$.  By the well-known comparison theorem, $\E|x(t)|^2\ge u(t)$ for $t \in [0,0.5\e]$, where
\begin{align} \label{7.9}
u(t) =  \bar z^2 +\int_0^t \big[ - (2+0.5 \bar z^2)u(s)  +u^2(s) \big] ds.
\end{align}
It is known (see, e.g., \cite{M97, MY06}) that equation (\ref{7.9}) has its explicit solution
$$
u(t) = \Big[ e^{(2+0.5\bar z^2)t} \Big( \frac{1}{\bar z^2} + \frac{1}{2+0.5\bar z^2} \big[  e^{-(2+0.5\bar z^2)t }- 1\big]   \Big) \Big]^{-1}
$$
on $t\in [0,0.5\e)$.  This implies $u(t)\to\8$ as $t\to 0.5\e$. Consequently, $\E|x(t)|^2$ and, hence $\E|x(t)|^4\to\8$ as $t\to 0.5\e$.  This is in contradiction to that $\E|x(t)|^4$ is finite and tends to 0 as $t\to\8$.  We can therefore conclude that
the controlled SDDE (\ref{7.4}) is NOT exponentially stable in the 4th moment no matter how small the $\e$ is, although its corresponding controlled SDE (\ref{7.2}) is exponentially stable in the 4th moment.

 \section*{Acknowledgements}
\par
Junhao Hu would like to thank the Natural Science Foundation of China (61876192) for the financial support.
\par
Wei Liu would like to thank the Natural Science Foundation of China (11701378, 11871343),
Chenguang Program supported by both Shanghai Education Development Foundation and Shanghai Municipal Education Commission (16CG50),
and Shanghai Gaofeng \& Gaoyuan Project for University Academic Program Development for their financial support.
\par
Feiqi Deng would like to thank the Natural Science Foundation of China (61503142, 61573156, 61733008) for their financial support.
\par
Xuerong Mao would like to thank
the Royal Society (WM160014, Royal Society Wolfson Research Merit Award),
the Royal Society and the Newton Fund (NA160317, Royal Society-Newton Advanced Fellowship),
and the EPSRC (EP/K503174/1) for their financial support.

\end{document}